\documentclass[11pt, reqno]{amsart}
\usepackage{amsfonts, amsmath, amssymb, latexsym, eucal, amscd} 
\begin{document}
 
%%%%%%%%%%%%%%%%%%% 

\numberwithin{equation}{section}

%%%%%%%%%%%%%%%%%%% 
\newtheorem{definition}[equation]{Definition}
\newtheorem{theorem}[equation]{Theorem}
\newtheorem{lemma}[equation]{Lemma}
\newtheorem{corollary}[equation]{Corollary}
\newtheorem{proposition}[equation]{Proposition}
\newtheorem{remark}[equation]{Remark}
\newtheorem{remarks}[equation]{Remarks}
\newtheorem{example}[equation]{Example}
\newtheorem{conjecture}[equation]{Conjecture}
\newtheorem{problem}[equation]{Problem}
\newtheorem{note}[equation]{Note}

\def\C{\mathbb C}
\def\K{\mathbb K} 
\def\Z{\mathbb Z}
\def\K{\mathbb F}
\def\N{\mathbb N}
\def\I{\mathbb I}
\def\fld{\mathbb K} 

\newcommand{\uqhat}{U_q(\widehat{\mathfrak{sl}}_2)}
\newcommand{\ot}{\otimes}

%%%%%
 \setcounter{secnumdepth}{2}

\title[]{\large The Universal Central Extension \\
of the Three-point $\mathfrak {sl}_2$ Loop Algebra}

\author[]{Georgia Benkart$^{\star}$} \address{Department of Mathematics \\ University
of Wisconsin \\  Madison, WI  53706, USA}
\email{benkart@math.wisc.edu}

\author[]{Paul Terwilliger \address{ \\ }}
\email{terwilli@math.wisc.edu}

\thanks{$^{\star}$Support from NSF grant \#{}DMS--0245082 is gratefully acknowledged.  \hfil \break 
Onsager Lie algebra.  
 }   
 
 \date{December 17, 2005}
 
 %\maketitle

 \begin{abstract} 
We consider the three-point loop algebra, 
$$L= \mathfrak {sl}_2\otimes \K \lbrack t, t^{-1}, (t-1)^{-1}\rbrack,$$
 where $\K$ denotes a field of characteristic $0$ and
 $t$ is
an indeterminate.  
The universal central extension 
$\widehat L$ of $L$ 
was determined by Bremner.   In this note,   we give
a presentation for
$\widehat L$
via generators and relations, which highlights a certain
symmetry over the alternating group $A_4$.  
To obtain our presentation of
$\widehat L$,  we use the realization of $L$ 
as the tetrahedron Lie algebra.
 \end{abstract} 
\medskip

 \maketitle
\section{
Introduction}

Throughout this paper
$\K$ will denote a field of characteristic 0.
Consider the Lie algebra
 $\mathfrak{sl}_2$ over $\K$, and   
for an indeterminate $t$, 
the polynomial algebra $\K[t]$ localized
at $t$ and $t-1$:  
\begin{eqnarray*}
\mathcal {A} = \K\lbrack t,t^{-1}, (t-1)^{-1}\rbrack.
\end{eqnarray*} 
%%%%%%%%%%%%%%%%%%%%%%%%%
%(Localization at any two
% $t-a_1$, $t-a_2$ for distinct $a_1,a_2 \in \K$ results in an  algebra
%isomorphic to $\mathcal A$,  so there is no loss of generality
%in assuming one value is $0$ and the other is $1$.)
%%%%%%%%%%%%%%%%%%%
% confusing distraction: since its not yet clear what we are trying to
%do, it seems premature to say there is no loss.
%%%%%%%%%%%%%%%%
 (Localization at any two
 $t-\alpha_1$, $t-\alpha_2$ for distinct $\alpha_1,\alpha_2 \in \K$ results in an  algebra
isomorphic to $\mathcal A$,  so there is no loss of generality
in assuming one value is $0$ and the other is $1$.)
The {\it loop algebra} corresponding to
$\mathfrak{sl}_2$
and $\mathcal A$ 
 is the Lie algebra 
% over $\K$ consisting of
%The loop algebra  $L$  is  the tensor product of 
%the Lie algebra $\mathfrak{sl}_2$ over $\K$ with $\mathcal A$, 
  \begin{eqnarray*}
L={\mathfrak{sl}_2}\otimes \mathcal {A}, 
\end{eqnarray*}
with product $[x \ot a, y \ot b] = [x,y] \ot ab$.

Our primary focus  here is on central extensions
of $L$, so we begin by recalling a few relevant definitions. 
A  {\it central extension} of a Lie algebra 
$\mathcal L$ is a 
pair $(\mathcal K,\pi)$ consisting of a Lie  algebra $\mathcal K$ and a surjective 
Lie algebra homomorphism $\pi:  \mathcal K  \rightarrow  {\mathcal L}$
whose kernel 
lies in the center of $\mathcal K$. 
Given central extensions
$(\mathcal K,\pi)$ and 
$(\mathcal K',\pi')$ of $\mathcal L$, 
by a {\it homomorphism} 
(resp. {\it isomorphism\/}) 
from 
$(\mathcal K,\pi)$
to $(\mathcal K',\pi')$ 
we mean 
%a central extension $\pi  : K \to 
%L$ to a central extension $\pi ' : K' \to L$ is a Lie algebra 
a homomorphism (resp. isomorphism)
of Lie algebras $\varphi : \mathcal K \to \mathcal K'$ such that
$\pi = 
\pi' \circ \varphi$.  

A central extension
$(\mathcal K, \pi)$ 
%$\widehat \pi : \widehat L \rightarrow  L$ is a {\it 
%universal central extension} if there exists a unique homomorphism 
of $\mathcal L$ is {\it universal}
whenever there exists a homomorphism
from 
$(\mathcal K, \pi)$  to any other central extension of $\mathcal L$.
A Lie algebra $\mathcal L$ possesses a universal central extension
if and only if 
$\mathcal L$ 
is perfect (i.e.  ${\mathcal L} = [{\mathcal L},{\mathcal L}]$),
and in this case, 
the 
universal central extension of $\mathcal L$ is unique up to  isomorphism.
It is well-established that  
the universal central extension plays a crucial
role in representation theory;
 one need only look at the examples of the affine, 
Virasoro, and toroidal Lie algebras for affirmation of this statement
(see \cite{Ka}, \cite{MP}).  

The loop algebra $L = \mathfrak {sl}_2 \ot \mathcal A$
is easily seen to be perfect; 
therefore,  $L$ has a universal central extension which we denote by 
$({\widehat L}, \pi)$.   Bremner \cite{Br} has given a detailed
description of 
$({\widehat L}, \pi)$.    He has shown that the
center of $\widehat L$ has dimension 2,   and he has given 
an explicit  basis and Lie bracket 
for 
${\widehat L}$.

%  He shows that the center of 
%${\widehat L}$ 
%coincides with
%the kernel of $\pi$  (this is always true for the u.c.e.)

Our goal here is to give a presentation for 
${\widehat L}$
via generators
%and relations that makes transparent its symmetry under
%the alternating group $A_4$. (..its symmetry..sounds like
%the full automorophism group)
and relations that highlights a certain symmetry over 
the alternating group $A_4$.  

Our point of departure is
the realization of $L$ as the tetrahedron
algebra  given by Hartwig and Terwilliger.
\medskip

\begin{definition}
\label{thm:ht}
\rm
\cite{HT}
Let $\boxtimes$ denote the Lie algebra over $\K$
defined by generators
$\{x_{i,j} \mid i,j \in \mathbb I, i \neq j\}$,  $\mathbb I = \{0,1,2,3\}$, 
and the following relations:
\begin{itemize}

\item[(i)]  For distinct $i,j \in \mathbb I$,
\begin{equation}\label{eq:1.1}x_{i,j} + x_{j,i} = 0. \end{equation}
\item[(ii)]  For mutually distinct $i,j, k \in \mathbb I$,
\begin{equation}\label{eq:1.2} [x_{i,j}, x_{j,k}] = 2 x_{i,j} + 2 x_{j,k}. \end{equation}
\item[ (iii)]  For mutually distinct $i,j,k,\ell \in \mathbb I$,
\begin{equation}\label{eq:1.3} [x_{i,j}[x_{i,j},[x_{i,j}, x_{k,\ell}]\,]\,]= 4 [x_{i,j},x_{k,\ell}].\end{equation}  
\end{itemize}
We call $\boxtimes$ the {\it tetrahedron algebra}. 
\end{definition}
\medskip

\begin{theorem}
\label{thm:boxt}
{\rm \cite{HT}}
The Lie algebras $\boxtimes$ and $L$ are isomorphic.
\end{theorem}

%Our strategy for describing a presentation of $\widehat L$
%is to use the  isomorphism  of Lie algebras
%$\sigma :\boxtimes \to L$ from \cite{HT}. (The sentence suggests
%that the use of sigma is the extent of the strategy)
%Then we will display a universal central extension (We wont know
%that it is universal until after we have shown hat sigma is an iso)
% $({\widehat \boxtimes}, \pi)$ of  
% $\boxtimes$,  with 
% ${\widehat \boxtimes}$ defined by
% generators and relations,  
% and a homomorphism
%of Lie algebras 
%${\widehat \sigma} 
%: {\widehat \boxtimes} \to {\widehat L}$
%such that the following diagram commutes:

We will obtain our  presentation of $\widehat L$
as follows.    First
we will display a central extension
$({\widehat \boxtimes}, \pi)$ of  
$\boxtimes$,  with 
 ${\widehat \boxtimes}$ defined by
 generators and relations.    Then we will modify 
the Lie algebra  isomorphism $\sigma :\boxtimes \to L$  given in \cite{HT}
 to obtain a homomorphism of Lie algebras, 
${\widehat \sigma} 
: {\widehat \boxtimes} \to {\widehat L}$, 
such that the following diagram commutes:
\begin{eqnarray}
\label{eq:cd}
%\[
\begin{CD}
{\widehat \boxtimes}  @>\pi>>
	  \boxtimes \\
	  @V{\widehat \sigma} VV          @VV\sigma V \\
	  {\widehat L}
	  @>\pi>>
		   L 
		    \end{CD}
%		\]
\end{eqnarray}
Using this and the universality of
$({\widehat L}, \pi)$,  we will argue
that 
${\widehat \sigma} 
: {\widehat \boxtimes} \to {\widehat L}$
is an isomorphism and thereby determine a presentation
of $\widehat L$ by generators and relations.
\medskip

\section{The isomorphism $\sigma :\boxtimes \to L$}

Let $\{e,f,h\}$ denote
the canonical basis for $\mathfrak {sl}_2$ having
product $[h,e] = 2e$, $[h,f] = -2f$, and $[e,f] = h$. 
The basis 
\begin{eqnarray*}
 X=2e-h, \qquad \qquad 
Y =-2f-h,
\qquad \qquad 
Z =h
\end{eqnarray*}
is more suitable for our purposes.  
Following
\cite{ITW} 
we call $X,Y,Z$ the {\it equitable basis}
for $\mathfrak{sl}_2$, since
%as the multiplication
%in terms of this basis is given by
 \begin{eqnarray*}% \label{eq:equit}%
\lbrack X,Y\rbrack = 2X+2Y,
\qquad \
\lbrack Y,Z\rbrack = 2Y+2Z,
\qquad \ \lbrack Z,X\rbrack = 2Z+2X.
\end{eqnarray*}

There exists a unique $\K$-algebra automorphism
$\prime $ of $\mathcal A$ such that
$t'=1-t^{-1}$.   This automorphism has order 3 
and satisfies
\begin{eqnarray}
&&t''=(1-t)^{-1},
\qquad \qquad t t'=t-1,  \label{eq:teq1}
\\
&&t' t''=t'-1, 
\qquad \qquad 
t'' t=t''-1,  \label{eq:teq2}
\end{eqnarray}  where $t'' = (t')'$.  
The relations in \eqref{eq:teq1} and \eqref{eq:teq2}  imply that the following is a basis for the $\K$-vector
space $\mathcal A$:
\begin{eqnarray*}
\lbrace 1 \rbrace \cup \lbrace t^i, (t')^i, (t'')^i \;|\;i \in \N\rbrace,
\end{eqnarray*} where $\N = \{1,2,3, \dots\}$.
\medskip

\begin{proposition}
\label{prop:ht}
{\rm{ (\cite[Thm. 11.5]{HT})}}  
There exists an isomorphism
of Lie algebras $\sigma : \boxtimes \to L$
that sends
\begin{eqnarray*}
&&x_{1,2}\;\to\; X\otimes 1,
\qquad \qquad x_{0,3}\;\to\;Y\otimes t + Z\otimes (t-1),
\\
&&x_{2,3}\;\to\;Y\otimes 1,
\qquad \qquad x_{0,1}\;\to\; Z\otimes t' + X\otimes (t'-1),
\\
&&x_{3,1}\;\to\;Z\otimes 1,
\qquad \qquad x_{0,2}\;\to\;X\otimes t'' + Y\otimes (t''-1)
\end{eqnarray*}
where $X,Y,Z$ is the equitable basis for 
$\mathfrak{sl}_2$.
\end{proposition}
\medskip

\section{A central extension of $\boxtimes$}

In this section,   we will construct a central extension $({\widehat \boxtimes}, \pi)$
of $\boxtimes$. Later we will show that
this extension is universal.   From now on,  
we identify the symmetric group $S_4$ with the group of
permutations of $\I$.  
\medskip
  
\begin{definition} 
\label{def:even}
\rm
Given a sequence $(i,j,k)$  of mutually distinct
elements of $\I$,   there exists a unique 
$\tau \in S_4$ such that $\tau(0)=i$,  
$\tau(1)=j$,
and 
$\tau(2)=k$.
The sequence  $(i,j,k)$ is said to be  {\it even} (resp. {\it odd})
whenever $\tau \in A_4$ (resp. $\tau \not\in A_4$),
where $A_4$ is the  alternating subgroup of $S_4$.
\end{definition} \medskip

\begin{definition} 
\label{def:shape}
\rm
A partition $p$ of $\I$ into two (disjoint) 
subsets, each with two elements, is said to have {\it shape} $(2,2)$.
\end{definition}

The set $P(\I)$ of all partitions of $\I$ of
shape $(2,2)$ has cardinality  3.  \medskip

\begin{definition}
\label{def:hatbox}
\rm
Let $\widehat \boxtimes $ denote
 the Lie
algebra over $\K$ defined by generators
\begin{eqnarray*}
\{X_{i,j} \mid i,j \in \I, i \neq j\} \cup 
\{C_p \mid p \in P(\I)\}
\end{eqnarray*}
and the following relations:
\begin{itemize}
\item[(i)]  For $p \in P(\I)$,
\begin{eqnarray*}
C_p \; {\mbox{ is  central}}.
\end{eqnarray*}
\item[(ii)] 
\begin{eqnarray*}
\sum_{p \in P(\I)}C_p=0.
\end{eqnarray*}
\item[(iii)]  For distinct $i,j \in \mathbb I$,
\begin{eqnarray*}
X_{i,j} + X_{j,i} = C_p, 
\end{eqnarray*}
where $p\in P(\I)$ consists of
$\lbrace i,j\rbrace$ and its complement in 
 $\I$.
\item[(iv)]  For mutually distinct $i,j, k \in \mathbb I$
such that $(i,j,k)$ is even,
\begin{eqnarray*}
[X_{i,j}, X_{j,k}] = 2 X_{i,j} + 2 X_{j,k}. 
\end{eqnarray*}
\item[ (v)]  For mutually distinct $i,j,k,\ell \in \mathbb I$,
\begin{eqnarray*}
[X_{i,j}[X_{i,j},[X_{i,j}, X_{k,\ell}]\,]\,]= 4 [X_{i,j},X_{k,\ell}].
\end{eqnarray*}  
\end{itemize}
\end{definition} \medskip

\begin{lemma}
\label{lem:boxtpi}
There exists a surjective 
homomorphism of Lie algebras
$\pi:  
{\widehat \boxtimes}
\to \boxtimes$ such that
\begin{eqnarray*}
\pi(X_{i,j})&=& x_{i,j} \qquad \qquad i,j \in \I,\; i\not=j,
\\
\pi(C_p)&=& 0 \qquad \qquad \quad  p \in P(\I).
\end{eqnarray*} 
\end{lemma}
\noindent {\it Proof:}
Compare the defining relations for
$\boxtimes$ and $\widehat \boxtimes$
given in Definitions
\ref{thm:ht} and
\ref{def:hatbox}. 
%That such a homomorphism $\pi$ exists is 
%readily seen  by simply comparing the defining relations for
%$\boxtimes$ and $\widehat \boxtimes$
%given in Definitions
%\ref{thm:ht} and
%\ref{def:hatbox}.   Uniqueness of the homomorphism  
%follows from the fact that the
%$X_{i,j}$ and the $C_p$ together generate
%$\widehat \boxtimes$.    The map $\pi$  is surjective, 
%since the $X_{i,j}$ map onto the 
%generators $x_{i,j}$ of  $\boxtimes$.
\hfill $\Box $  \medskip

\begin{lemma}
\label{lem:detail}  
For mutually distinct 
$i,j,k \in \I$ such that $(i,j,k)$ is odd,
%in the
%sense of
%Definition 
%\ref{def:even}.  Then in the algebra $\widehat \boxtimes$ we have 
 in the algebra $\widehat \boxtimes$ we have 
\begin{eqnarray}
\lbrack X_{i,j}, X_{j,k}\rbrack
=2X_{i,j}+2X_{j,k}+2C_p,
\label{eq:xxc}
\end{eqnarray}
where $p \in P(\I)$ consists of
$\lbrace i,k\rbrace$ and its complement 
in $\I$.
\end{lemma}
\noindent {\it Proof:} The sequence
$(k,j,i)$ is even since
$(i,j,k)$ is odd.  Therefore, by Definition
\ref{def:hatbox}\,(iv),
\begin{eqnarray*}
\lbrack X_{k,j}, X_{j,i}\rbrack
=2X_{k,j}+2X_{j,i}.
\label{eq:xxnoc}
\end{eqnarray*}
Evaluating this using (i)-(iii) of
 Definition
\ref{def:hatbox}, 
we obtain 
(\ref{eq:xxc}).
\hfill $\Box $ \medskip

\begin{lemma}
\label{lem:coincide}
The following subspaces of $\widehat \boxtimes$ coincide:
\begin{itemize}
\item[{\rm (i)}] 
the kernel of $\pi$,
\item[\rm{(ii)}] 
$\mbox{\rm Span}\{C_p \mid p \in P(\I)\}$,
\item[\rm{(iii)}] 
the center of $\widehat \boxtimes$.
\end{itemize}
\end{lemma}
\noindent {\it Proof:} 
Set
${\mathcal C}=\mbox{Span}\{C_p \mid p \in P(\I)\}$.
We first show
${\mathcal C}=\mbox{Ker}(\pi)$.
We have 
 ${\mathcal C}\subseteq \mbox{Ker}(\pi)$ by
Lemma
\ref{lem:boxtpi}. To establish equality,
observe that $\mathcal C$ is an ideal
in $\widehat \boxtimes$, 
and let $\pi': {\widehat \boxtimes} \to {\widehat \boxtimes}/{\mathcal C}$
denote canonical surjection with kernel $\mathcal C$. 
Since $\pi'(C_p)= 0$ for $p\in P(\I)$,  it follows from  Definition
\ref{def:hatbox} and Lemma
\ref{lem:detail} that  the elements
$\{\pi'(X_{i,j}) \mid i,j \in \mathbb I,\; i \neq j\}$
satisfy the defining relations 
(\ref{eq:1.1})--(\ref{eq:1.3})
for $\boxtimes$. Therefore,  there exists
a Lie algebra homomorphism
$\gamma :\boxtimes \to {\widehat \boxtimes}/{\mathcal C}$ such that
$\gamma(x_{i,j})= \pi'(X_{i,j})$
for all distinct $i,j \in \I$. From the construction, 
the following diagram commutes:
\[
\begin{CD}
{\widehat \boxtimes}  @>\pi>>
	  \boxtimes \\
	  @V{\rm {id}} VV          @VV{\gamma}V \\
	  {\widehat \boxtimes}
	  @>\pi'>>
		   {\widehat \boxtimes}/{\mathcal C} 
		    \end{CD}
		    \]
We may now argue
\begin{eqnarray*}
{\mathcal C} &=& \mbox{Ker}(\pi')
\\
 &=& \mbox{Ker}(\gamma \circ \pi)
\\
&\supseteq & 
 \mbox{Ker}(\pi), 
\end{eqnarray*}
which implies that 
${\mathcal C}=
 \mbox{Ker}(\pi)$.
Next we prove that the center
$Z(\widehat \boxtimes)={\mathcal C}$.
We have  
${\mathcal C}\subseteq Z(\widehat \boxtimes)$
by Definition
\ref{def:hatbox}\,(i). To obtain the
reverse inclusion, it suffices to  
show that the image of 
$Z(\widehat \boxtimes)$ under $\pi$ is zero.
This image is contained in  
$Z(\boxtimes)$ since
$\pi$ is surjective.
But $\boxtimes $ is
isomorphic to $L$ and
$Z(L)=0$,  
so
$Z(\boxtimes)=0$. 
Therefore the image of 
$Z(\widehat \boxtimes)$ under $\pi$ is zero and
consequently 
$
Z(\widehat \boxtimes)
\subseteq
{\mathcal C}$.
From these comments,  we find that
$Z(\widehat \boxtimes)={\mathcal C}$.
%The map $\pi: {\widehat \boxtimes} \to \boxtimes$
%is a surjection with kernel $J$,  and 
%the image of $Z(\widehat \boxtimes)$
%under $\pi$ is contained in
%$Z(\boxtimes)$.   But $\boxtimes $ is
%isomorphic to $L$ and
%$Z(L)=0$, 
%so
%$Z(\boxtimes)=0$.   Thus, $Z(\widehat \boxtimes) \subseteq J$,
%and since the reverse  inclusion holds as well, 
%$J= Z(\widehat \boxtimes)$.
\hfill $\Box $ \medskip

\begin{corollary}
\label{cor:ce}
The pair $({\widehat \boxtimes}, \pi)$ is a central extension
of $\boxtimes$.
\end{corollary}
\noindent {\it Proof:} 
The map $\pi: {\widehat \boxtimes} \to \boxtimes$
is a surjective 
homomorphism of Lie algebras whose kernel is
contained in the center of $\widehat \boxtimes$.
\hfill $\Box $

%\noindent It is so named because we may think of a tetrahedron with vertices
%indexed by the elements of $\mathbb I$ and identify the generator $X_{i,j}$ 
%with an edge from
%$i$ to $j$ for $i \neq j$. 
%In \cite{HT},   it is shown that the tetrahedron algebra  is a direct sum of three copies of
%the Onsager Lie algebra.   The Onsager algebra was first
%introduced in a seminal paper  \cite{O}  in which the free energy of the
%two-dimensional Ising model was computed exactly.  
%It has been investigated extensively in connection with solvable
%lattice models in physics, and on the mathematical side, in connection
%with representation theory, Kac-Moody Lie algebras, tridiagonal pairs,
%and partially orthogonal polynomials.    In \cite{P}, Perk showed that
%the Onsager Lie algebra $\mathcal O$ has a simple presentation by
%generators $A,B$ which satisfy the relations
%
%\begin{eqnarray*} [A,[A,[A,B]\,]\,] &=& 4[A,B] \\  
%{[B,[B,[B,A]\,]\,]} &=& 4[B,A]. \end{eqnarray*}
%Hence,  each pair $X_{i,j}$, $X_{k,\ell}$  with distinct subscripts $i,j,k,\ell$ from $\mathbb I$
%determines an Onsager subalgebra of the loop algebra $L$.   
%

%The realization of $L$
%the loop algebra $L = \mathfrak{sl}_2 \ot \mathcal A$ as
%as the tetrahedron algebra makes its symmetry under the symmetric group $S_4$
%transparent.  It seems natural to ask for
%a similar realization of $\widehat L$, and this
%is what we will provide later in the paper.
%%%%%%%
%the universal central extension of $L$
%that makes evident its symmetry and exhibits its connections with the Onsager algebra. 
%That is what we aim to accomplish here.  
\medskip

\section{
 The universal central extension of
$L =   
\mathfrak {sl}_2 \ot \mathcal A$}

In this section,  we give a detailed
description of the universal central extension
$({\widehat L}, \pi)$ of $L$. In the next section
we will use this
description  to define the map
${\widehat \sigma} : 
{\widehat \boxtimes}
\to {\widehat L}$
mentioned in Section 1.

\medskip
 
 By results of Kassel \cite{K} (see also
\cite{KL} and \cite{BK}),  the universal central extension of  $\mathcal L: = \mathfrak g \ot \mathcal B$
for any finite-dimensional complex simple Lie algebra $\mathfrak g$  and any
commutative, associative algebra $\mathcal B$ with 1 is  obtained from 
$\mathcal L$ 
 by adjoining the  K\"ahler differentials
modulo the exact forms of $\mathcal B$.  
This description 
          %of $\widehat {\mathcal L}$
enabled Bremner \cite{Br}  to show that the universal central extension of any $n$-point
loop algebra over $\mathfrak g$  has an  
$(n-1)$-dimensional kernel.
The same argument  works
over any field $\K$ of characteristic 0.    Applying this result to our loop algebra    
$L = \mathfrak {sl}_2 \ot \mathcal A$, we see that $\dim_\K  \widehat L/L = 2$.  

An alternative  description of $\widehat L$ 
can be found in [ABG].
Let $\mathcal S$ denote the subspace of $\mathcal A \ot \mathcal A$ 
%generated
spanned by the elements $a \ot b + b \ot a$ and $ab \ot c + bc \ot a + ca \ot b$ for all $a,b,c \in
\mathcal A$.   Let $\langle \mathcal A, \mathcal A\rangle  =
(\mathcal A \ot \mathcal A)/ \mathcal S$, and write
$\langle a, b \rangle = (a \ot b) + \mathcal S$.  From the construction
we know that 
\begin{eqnarray}&& 
\langle a,b\rangle + \langle b,a \rangle = 0, \label{eq:ucerel1}\\
&&\langle ab, c \rangle + \langle bc, a \rangle + \langle ca, b \rangle = 0 \label{eq:ucerel2}
\end{eqnarray}
for all $a,b,c \in \mathcal A$.  Then  
$$\widehat L = \big(\mathfrak {sl}_2 \ot \mathcal A\big) \oplus
\langle \mathcal A, \mathcal A\rangle,$$ 
where 
$\langle \mathcal A, \mathcal A\rangle$ is central and
$$\lbrack x\otimes a, y\otimes b \rbrack  
= \lbrack x,y\rbrack \otimes ab + (x\,|\,y)\langle a,b\rangle
$$
for all $x,y \in 
\mathfrak {sl}_2$
and $a,b\in \mathcal A$.
%is a Lie algebra, and $\widehat L$  is, in fact, the universal central extension of $L$ (see \cite{ABG}). 
Here $(x\,|\,y)$ denotes the Killing form of
$\mathfrak {sl}_2$.  Thus, finding $\widehat L$ amounts to computing $\langle \mathcal A, \mathcal A
\rangle$ explicitly.

Using the relations 
\eqref{eq:ucerel1}, \eqref{eq:ucerel2},   it is not difficult to show,  just
as in the affine case (see \cite{Ka},  for example), that

\begin{equation} \langle f^m, f^n \rangle  = m \delta_{m+n,0} \langle f, f^{-1} \rangle \end{equation}
for $f = t,t'$,  or $t''$ and all integers $m,n$, where $\delta$ is the Kronecker delta. 
Letting $g =f'$ and using
$f'=1 - f^{-1}$, we have  
\begin{eqnarray}   \langle f^m, g^n \rangle  
%&=& \sum_{k=0}^n (-1)^k  {n \choose k} \langle f^m, f^{-k} \rangle \\
&=& \sum_{k=0}^n (-1)^k  \binom{n}{k} \langle f^m, f^{-k} \rangle \\
%&=& m (-1)^m {n \choose m} \langle f, f^{-1} \rangle
&=& m (-1)^m \binom{n}{m} \langle f, f^{-1} \rangle
\nonumber \end{eqnarray}
for nonnegative integers $m,n$.
%As  $t'  = 1- t^{-1}$, $t'' = 1- (t')^{-1} $, and $t'  = 1- (t'')^{-1}$,
%we see that 
%$\langle \, , \, \rangle$ is completely determined by these relations.  
%Moreover,   since $t t' = t-1 = (t'')^{-1}$,
%\ $t' t'' = t'-1 = t^{-1}$, and 
%$t'' t = t'' -1 = (t')^{-1}$,  we have
From (\ref{eq:teq2}) we find that 
$tt't''=-1$; using this and
(\ref{eq:ucerel2})
we obtain
\begin{eqnarray*} \langle t'', (t'')^{-1} \rangle 
&=&  -\langle t'', tt' \rangle \\
&=&   \langle t, t' t'' \rangle + \langle t', t''t \rangle\\
&=&  -\langle t,t^{-1} \rangle - \langle t', (t')^{-1} \rangle. 
\end{eqnarray*}
It follows from these computations
  that  
%$\widehat L$ is spanned modulo $L$
$\langle \mathcal A, \mathcal A\rangle$  is spanned by 
$\langle t,t^{-1} \rangle$ and $\langle t',(t')^{-1} \rangle$
and  that 
\begin{eqnarray*}
\langle t,t^{-1} \rangle +\langle t',(t')^{-1} \rangle + 
\langle t'', (t'')^{-1} \rangle = 0.
\end{eqnarray*}
 Since
 ${\widehat L}/L$ has dimension 2,  the space
$\langle \mathcal A, \mathcal A\rangle$ has dimension 2.
Consequently 
$\langle t,t^{-1} \rangle$ and $\langle t',(t')^{-1} \rangle$
form a basis for 
$\langle \mathcal A, \mathcal A\rangle$.

\medskip
 
The next result is now apparent.
\medskip

\begin{theorem}
\label{thm:brem}
{\rm { (\cite{Br},
\cite{ABG})}}
Let $C$ denote a two-dimensional vector space over $\K$.
Let 
$\mathfrak c,
 \mathfrak c'
$ 
denote a basis for $C$ and define
 $\mathfrak c''$
 so that
\begin{eqnarray*}
 \mathfrak c
+ 
 \mathfrak c'
+
\mathfrak c'' = 0.
\end{eqnarray*}
Then the following $\rm{(i)}$--$\,\rm{(iii)}$ hold.
\begin{enumerate}
\item[{\rm (i)}] 
There exists a Lie algebra  
  \begin{eqnarray*}
{\widehat L} = L \oplus C
\end{eqnarray*}
with product
\begin{eqnarray*}
&&\lbrack {\widehat L}, C\rbrack=0,
\\
&&\lbrack x\otimes a, y\otimes b \rbrack  
= \lbrack x,y\rbrack \otimes ab + (x\,|\,y)\langle a,b\rangle
\end{eqnarray*}
for $x,y \in 
 \mathfrak {sl}_2$
and $a,b\in \mathcal A$, where 
 $(x\,|\,y)$ is the Killing form for
$\mathfrak {sl}_2$,  and where 
$\langle \,,\,\rangle : \mathcal A \times \mathcal A \to C$ is $\K$-bilinear
and satisfies 

\bigskip \smallskip
\centerline{
\begin{tabular}[t]{c|cccc}
        $\langle\, , \,\rangle$  
       & 1  
       & $t^j$  
         & $(t')^j$  
         & $(t'')^j$  
	\\
          \\
	\hline 
\\
$1$ & $0$ & $0$ & $0$ &$0$
\\
\\
$t^i$
&
$0$
& 
$0$
&
$(-1)^i i 
\binom{j}{i} 
\mathfrak c $
&
$(-1)^{j+1} j 
\binom{i}{j} 
\mathfrak c'' $
\\
\\
$(t')^i$
& 
$0$
&
$(-1)^{j+1} j 
\binom{i}{j} 
\mathfrak c$
&
$0$
&
$(-1)^i i 
\binom{j}{i} 
\mathfrak c' $
\\
\\
$(t'')^i$
& 
$0$
&
$(-1)^i i 
\binom{j}{i} 
\mathfrak c'' $
&
$(-1)^{j+1} j 
\binom{i}{j} 
\mathfrak c'$
&
$0$
\end{tabular}}

\bigskip \smallskip
\noindent for $i,j\in \N$.
\smallskip 

\item[{\rm (ii)}] 
There exists a homomorphism of Lie algebras
$\pi : {\widehat L}\to L$ that has kernel $C$ and
acts as the identity on $L$.  \smallskip

\item[{\rm (iii)}] 
The pair 
$({\widehat L}, \pi)$ is the universal central extension of
$L$. 
\end{enumerate}
\end{theorem}
\medskip

\section{A homomorphism ${\widehat \sigma}: {\widehat \boxtimes} \to
{\widehat L}$}

\begin{lemma}
\label{thm:widesig}
There exists a unique Lie algebra homomorphism
${\widehat \sigma}: {\widehat \boxtimes} \to
{\widehat L}$ specified  by the following tables.

\bigskip
\centerline{
\begin{tabular}[t]{c|c}
        $p$  
       & {\rm Image of $C_p$ under 
${\widehat \sigma}$}
	\\
	\hline 
$\lbrace 0,1\rbrace \lbrace 2,3 \rbrace$ & 
$-4\mathfrak{c}^{\prime \prime}$
\\
$\lbrace 0,2\rbrace \lbrace 1,3 \rbrace$ & 
$-4\mathfrak{c}$
\\
$\lbrace 0,3\rbrace \lbrace 1,2 \rbrace$ & 
$-4\mathfrak{c}^{\prime}$
\end{tabular}}

\bigskip
\centerline{
\begin{tabular}[t]{c|c|c}
        $i$ & $j$  
       & {\rm Image of $X_{i,j}$ under   
${\widehat \sigma}$}
\\
	\hline 
$1$ &$2$ & 
$X\otimes 1 -4\mathfrak{c}'$
\\
$2$& $1$ & 
$-X\otimes 1$
\\
$2$&$3$ & 
$Y\otimes 1 -4\mathfrak{c}''$
\\
$3$& $2$ & 
$-Y\otimes 1$
\\
$3$&$1$ & 
$Z\otimes 1 -4\mathfrak{c}$
\\
$1$&$3$ & 
$-Z\otimes 1$
\\
$0$&$3$ & 
$Y\otimes t+Z\otimes (t-1) +4\mathfrak{c}$
\\
$3$&$0$ & 
$-Y\otimes t-Z\otimes (t-1) +4\mathfrak{c}''$
\\
$0$&$1$ & 
$Z\otimes t'+X\otimes (t'-1) +4\mathfrak{c}'$
\\
$1$&$0$ & 
$-Z\otimes t'-X\otimes (t'-1) +4\mathfrak{c}$
\\
$0$&$2$ & 
$X\otimes t''+Y\otimes (t''-1) +4\mathfrak{c}''$
\\
$2$&$0$ & 
$-X\otimes t''-Y\otimes (t''-1) +4\mathfrak{c}'$
\end{tabular}}, 

\bigskip
\noindent 
where $X,Y,Z$ is the equitable basis
for 
${\mathfrak sl}_2$.
\end{lemma}
\noindent {\it Proof:}
It is routine to verify that the elements
in the above tables satisfy the defining
relations
for
$\widehat \boxtimes$ given in 
Definition
\ref{def:hatbox}.
\hfill $\Box $  \medskip

\begin{lemma}
\label{lem:comd} The diagram in 
(\ref{eq:cd})
commutes.
\end{lemma}
\noindent {\it Proof:}
By Definition
\ref{def:hatbox}, 
the following is a generating set for
${\widehat \boxtimes}$:
\begin{eqnarray*}
\{X_{i,j} \mid i,j \in \I, i \neq j\} \cup 
\{C_p \mid p \in P(\I)\}.
\end{eqnarray*}
We chase these generators around the diagram
using 
the maps  in
Proposition \ref{prop:ht},
Lemma \ref{lem:boxtpi},
Theorem \ref{thm:brem}(ii),
and Lemma
\ref{thm:widesig}.
For each generator the image under the composition
${\pi} \circ {\widehat \sigma}$
coincides with the image under
the composition
${\sigma} \circ {\pi}$.
\hfill $\Box $  \medskip

\begin{theorem}
\label{thm:iso}
The Lie algebra homomorphism
${\widehat \sigma}: {\widehat \boxtimes} \to
{\widehat L}$ from
Lemma
\ref{thm:widesig}  is an isomorphism.
\end{theorem}
\noindent {\it Proof:}
We first show that $ {\widehat \sigma}$ is surjective. 
The map
$\pi : 
{\widehat \boxtimes} \to
{\boxtimes}$
is surjective,  and the 
map
$\sigma : \boxtimes \to L$
is an isomorphism,  so
the composite map
$\sigma \circ \pi:
{\widehat \boxtimes} \to
L$ is surjective.
Therefore  by Lemma
\ref{lem:comd},  the composition
$\pi \circ {\widehat \sigma} :
{\widehat \boxtimes} \to
L$ is  surjective.
The kernel of $\pi : {\widehat L} \to L$
is the space $C$ from
Theorem
\ref{thm:brem}, and $C$ is contained in
the image of
$ {\widehat \sigma}$, so 
$ {\widehat \sigma}$ is surjective. 
We now argue  that
$ {\widehat \sigma}$
is injective.
As before, set  $\mathcal C=
\mbox{Span}\{C_p \mid p \in P(\I)\}$.
The map
$\pi : 
{\widehat \boxtimes} \to
{\boxtimes}$ has kernel 
$\mathcal C$
 by  Lemma
\ref{lem:coincide},   
and the map $\sigma :\boxtimes \to L$
is an isomorphism, so  
the composition  
$\sigma \circ \pi:
{\widehat \boxtimes} \to
L$ has kernel  $\mathcal C$.
By this and  Lemma
\ref{lem:comd},  the composition $\pi \circ {\widehat \sigma} :
{\widehat \boxtimes} \to L$
has kernel  
$\mathcal C$.
Consequently,   the kernel of
$\widehat \sigma$ is contained in
$\mathcal C$.
From the first table of
Lemma
\ref{thm:widesig},
and the fact that ${\mathfrak c}, 
{\mathfrak c'}$ form a basis for $C$,
we find 
that the restriction of
$\widehat \sigma$ to 
$\mathcal C$ is injective.
Therefore 
$\widehat \sigma$ is injective and hence
an 
 isomorphism.   
\hfill $\Box $

\medskip 

\begin{corollary}
The central extension 
$({\widehat \boxtimes}, \pi)$ of $\boxtimes$ 
given in
Corollary
\ref{cor:ce}
is 
 universal.
\end{corollary}
\noindent {\it Proof:}
By Lemma 
\ref{lem:comd},  the diagram in 
(\ref{eq:cd})
commutes.
The map 
$\sigma $ is an isomorphism
by Proposition
\ref{prop:ht}, 
and
$\widehat \sigma$ is an isomorphism
by
Theorem
\ref{thm:iso}. The result follows
from this,   since 
$({\widehat L}, \pi)$ is the universal central
extension of $L$.
\hfill $\Box $
\medskip

\section{The defining relations for $\widehat \boxtimes$ revisited}

Definition 
\ref{def:hatbox} gives the defining relations 
 for the Lie algebra $\widehat \boxtimes$.
In this section,  we re-express these relations, this time
using a notation that makes explicit the role of the alternating group $A_4$.

 We will view  $S_4$ as acting on $\I$  from the right; this means
that when applying  a product $\alpha \beta $,  we first
apply
$\alpha $ and then $\beta$.
We consider the following normal subgroup of $S_4$:
\begin{eqnarray*}
N &=& \lbrace  (01)(23), \; (02)(31), \;(03)(12), \;e\rbrace. 
\end{eqnarray*} 
This subgroup is contained in $A_4$ and is therefore
 a normal subgroup of $A_4$.
 Let $N'$ denote the set of nonidentity elements of $N$.
For each $\eta \in N'$ let $\lbrack \eta \rbrack$
denote the partition of $\I$ consisting of the orbits
of $\eta$. Note that the map $\eta \to \lbrack \eta \rbrack$
is a bijection from $N'$ to $P(\I)$.
Two elements of $A_4$,
\begin{eqnarray*}
\zeta = (01)(23),
\qquad \qquad \vartheta=(012),
\end{eqnarray*}

\noindent play a distinguished role in our computations.  The first belongs to $N'$,
while the second one does not.   Together they generate $A_4$.    \medskip

\begin{theorem}
$\widehat \boxtimes$ is isomorphic to
the Lie algebra over $\K$ that has generators
\begin{eqnarray*}
X_\alpha, \;  \ C_\eta \qquad \qquad \alpha \in A_4, \; \ \eta  \in N' 
\end{eqnarray*}
and the following relations:
\begin{enumerate}
\item[{\rm (i)}] For $\eta \in N'$,
\begin{eqnarray*}
C_\eta \; {\mbox{ is  central}}.
\end{eqnarray*}
\item[\rm{(ii)}] 
\begin{eqnarray*}
\sum_{\eta \in N'}C_\eta=0.
\end{eqnarray*}
\item[\rm{(iii)}] For $\alpha \in A_4$,
\begin{eqnarray*}
X_\alpha +X_{\zeta \alpha} = C_{\alpha^{-1}\zeta \alpha}.
\end{eqnarray*}
\item[\rm{(iv)}] For $\alpha\in A_4$,
\begin{eqnarray*}
\lbrack X_\alpha, X_{\vartheta \alpha}\rbrack  = 
2X_\alpha +2X_{\vartheta \alpha}.
\end{eqnarray*}
\item[{\rm (v)}] For $\alpha \in A_4$ and for $\eta \in N'$, \ $\eta \not=\zeta$,
\begin{eqnarray*}
\lbrack X_{\alpha},
\lbrack X_{\alpha},
\lbrack X_{\alpha},
X_{\eta \alpha}\rbrack \rbrack \rbrack= 
4 \lbrack X_{\alpha},
X_{\eta \alpha}\rbrack.
\end{eqnarray*}
\end{enumerate}
An isomorphism with the presentation 
in Definition
\ref{def:hatbox} is given by
\begin{eqnarray*}
X_\alpha &\to& X_{\alpha(0),\alpha(1)}   \qquad \alpha \in A_4
\\
C_\eta &\to& C_{\lbrack \eta \rbrack} \qquad     \qquad \eta \in N'.  
\end{eqnarray*}
\end{theorem}
\noindent {\it Proof:} Up to notation, 
the above presentation of 
$\widehat \boxtimes$ is the same
as the one
given in Definition
\ref{def:hatbox}.
\hfill $\Box $
\medskip

\section {Concluding Remarks}
 
%In \cite{HT},   it is shown that the tetrahedron algebra  is
%a direct sum of three subalgebras,
%each of which 
%is isomorphic to the Onsager Lie algebra.
We conclude with some comments relating our results to 
the Onsager Lie algebra.
This algebra was introduced
in a seminal paper  \cite{O}  in which the free energy of the
two-dimensional Ising model was computed exactly.  
Since then it has been widely investigated 
by both the physics and mathematics communities 
in connection with solvable
lattice models, 
representation theory, Kac-Moody Lie algebras, tridiagonal pairs,
and partially orthogonal polynomials. 
In \cite{P}, Perk showed that
the Onsager Lie algebra  has a presentation by
generators $A,B$ and the following relations:
\begin{eqnarray*} [A,[A,[A,B]\,]\,] &=& 4[A,B] \\  
{[B,[B,[B,A]\,]\,]} &=& 4[B,A]. \end{eqnarray*}
Let $\Omega$ (resp. $\Omega'$) (resp. $\Omega''$)
denote the subalgebra of $\boxtimes$
generated by $x_{0,1}$ and $x_{2,3}$
(resp. $x_{0,2}$ and  $x_{1,3}$)
(resp. $x_{0,3}$ and $x_{1,2}$).
It was shown in 
\cite{HT} that each of the Lie algebras
$\Omega$, $\Omega'$, $\Omega''$ is isomorphic
to the Onsager algebra, and that
$\boxtimes$ is their direct sum.
Our results lead to a similar decomposition of 
$\widehat \boxtimes$ as follows:
Let
$\mathcal O$ 
(resp. ${\mathcal O}'$) 
(resp. ${\mathcal O}''$)
denote the subalgebra of
$\widehat \boxtimes$
generated by $X_{0,1}$ and $X_{2,3}$
(resp.
 $X_{0,2}$ and $X_{1,3}$)
(resp. $X_{0,3}$ and $X_{1,2}$). 
Using the commuting 
diagram
(\ref{eq:cd}) 
and the tables in Lemma \ref{thm:widesig}, 
we find that each of the algebras 
$\mathcal O$,
${\mathcal O}'$, 
${\mathcal O}''$ is isomorphic to the
Onsager algebra,
and that
$\widehat \boxtimes$ is the direct sum
$\mathcal O+
{\mathcal O}'+
{\mathcal O}''+{\mathcal C}$,
where
$\mathcal C := \mbox{Span}\{C_p \mid p \in P(\I)\}$ is
the two-dimensional center of $\widehat \boxtimes$.
Using the isomorphism  ${\widehat \sigma}:{\widehat \boxtimes}
\to {\widehat L}$ in Theorem \ref{thm:iso}, 
we obtain a corresponding decomposition for the universal
central extension $\widehat L$ 
of the loop algebra $L$.

\end{document}